\documentclass{article}
\usepackage{stmaryrd}
\usepackage{mathrsfs}
\usepackage{amsmath}
\usepackage{amssymb,amsfonts}
\usepackage[dvips]{graphicx}
\def\go{\rightarrow}
\def\half{{1\over 2}}

\begin{document}

\title{The Concentration and Stability of the Community Detecting Functions on Random Networks}
\author{Weituo Zhang, Chjan C. Lim}
%\affil{Mathematical Sciences, RPI, 110 8th St., Troy, NY, 12180}
%\affiliation{Mathematical Sciences, RPI, 110 8th St., Troy, NY, 12180. email: limc@rpi.edu}

\maketitle
\begin{abstract}
We propose a general form of community detecting functions for finding the communities or the optimal partition of a random network, and examine the concentration and stability of the function values using the bounded difference martingale method. We derive LDP inequalities for both the general case and several specific community detecting functions: modularity, graph bipartitioning and q-Potts community structure. We also discuss the concentration and stability of community detecting functions on different types of random networks: the sparse and non-sparse networks and some examples such as ER and CL networks.
\end{abstract}

\setlength{\parskip}{2ex}

\section{Introduction}
One of the main problems in network study is finding the communities or the optimal partition for a given network. A standard approach is to design a community detecting function of network partitions which achieves its extremum when the partition is optimal, or conversely the optimality of the partition is defined by a reasonably designed community detecting function.

There are several important applications of this approach. One is to find the community structures of social networks. A well known community detecting function for this application is the Modularity \cite{7} proposed by M.E.J. Newman. Another application is called circuit partioning in designing a computer system \cite{8}, in which circuits must be partitioned into groups so that the number of signals crossing the partition boundaries is minimized. The community detecting function used in this application is a variation of the function used in the graph bipartitioning problem \cite{5}. In this paper, we propose a general form of community detecting functions for which the functions mentioned above and many others are considered as specific cases.

Previous studies show that optimizing the community detecting function is in general an NP-complete problem \cite{9}. Algorithms such as simulated annealing is developed and tested on specific networks\cite{8}. However, in theoretical studies, it is nearly impossible to obtain an exact solution for a single network. What we usually have is the average value of the community detecting function over a random network ensemble\cite{5,6}. This gap between the computational and theoretical points of views arise the necessity to study the concentration and stability of the community detecting functions.

Roughly speaking, the stability of the community detecting function means that the change of the community detecting function is small when the corresponding network structure undergoes a small perturbation; and the concentration means that the theoretically predicted average value of the community detecting function becomes more precise when the system size grows larger.
In this paper, we derive LDP type inequalities to illustrate the both aspects: the fluctuation of the community detecting function in a fixed network ensemble and its asymptotic behavior when the system size goes infinity.

The concentration and stability of a community detecting function are  important from the following point of view:

Firstly, some times it is helpful to consider the given network as a sample randomly picked from a designed network ensemble.  The ensemble should be easy to analyze, catch some important features of the given network such as the average degree or the degree distribution, but neglect the detailed structural informations by randomization. In this scenario, the concentration and stability estimate the departure of the function value of the given network from the ensemble average.

Secondly for real world networks such as the social network and the Internet, we have problems like lack of information, uncertainty of the environment, and changing of the network by time. The optimization problem on this kind of networks is by the nature on an network ensemble. Since this ensemble is usually hard to analyze, simulation results of specific networks is in turn taken to estimate the ensemble average. In this scenario, the concentration and stability come from the  problem itself. For example, in the circuit partioning, if we have already found an optimal circuit partition under a given signal flow configuration, the concentration and stability conclude that this partition is still ``good'' if is not optimal under similar configurations.

Finally the concentration and stability give a measure for how a specific network deviates from an ensemble. For example, we calculate the modularities for a given network $G$ and the ER network ensemble. These two values by themselves do not tell you whether the communities are well defined in $G$. However, the LDP inequality from the concentration gives a bound of the probability that a network randomly picked from the ER network ensemble has smaller modularity than $G$ by chance. If the probability is small enough, it is statistically significant that $G$ is not picked from the ER ensemble considering the community structure.

Although not essential to the analysis in this paper, the general form of the community detecting function can be physically interpreted as a Hamiltonian of a spin system. We study the concentration and stability through a classical approach invented by J. Spencer et al \cite{Shamir1987}: consider an exploring process (edge by edge or node by node) of the network, construct a Doob martingale and take advantages of martingale inequalities.

We apply our results to several special community detecting functions: modularity, graph bipartitioning and q-potts community structure. We also discuss the problems on different classes of networks such as Erdos-Renyi (ER) and Chung-Lu (CL) networks.
Considering the asymptotic behavior of the number of edges when the network size grows, we also classify the network models into the sparse and non-sparse networks. We only study the ``sparse'' network defined by a constant upper bound of the degree, under which we proved a very general concentration result for all these problems. For the non-sparse network we derive our results only for specific cases.

\section{Background}
The underlying martingale inequality (Azuma's inequality) can be tracked back to Chebyshev inequality.  Bernstein inequality \cite{Bernstein1924,Bernstein1937} named after S.N.Bernstein is considered as a modification of Chebyshev inequality and gives a exponential decreasing probability upper bound. Then Hoeffding (1963) \cite{Hoeffding1963} invents and proves the first version of so called ``Hoeffding-Chernoff bound'' which gives a very general probability upper bound for sum of i.i.d.'s. In the same paper, Hoeffding also propose a slightly different version to the case where the random variables are not necessarily identical but uniformly bounded, and this version is usually called "Hoeffding's inequality". Azuma \cite{Azuma} extends the independent variables in Hoeffding's inequality to martingale differences and obtain the "Azuma's inequality". It is a great improvement for us, because the independence between the attendances of edges is usually unavailable in the network study.

W.T.Rhee and M. Talagrand (1987) \cite{Rhee1987} apply Azuma's inequality to the NP hard optimization problems. They give stochastic models for the Bin Packing problem and Traveling-Salesman problem as two examples, and assume a sequence of optimal solutions of a growing system to construct the martingale.

Shamir, Spencer (1987) \cite{Shamir1987} use Azuma's inequality on the problem of chromatic number of random graph. C. McDiarmid's study (1989) \cite{McDiarmid1989} summarizes this technique on random graph as the method of bounded differences and discusses some extensions such as the isoperimetric inequalities. Then several statistics of random graphs are studied such as the average distance \cite{Chung2002a}, connected component size \cite{Chung2002b} and number of triangles \cite{Kim2004}. More generalized martingale inequalities \cite{Godbole1998,Kim2004} and  models of random graphs \cite{Chung2002a,Chung2002b} are considered in these studies.

Our study combines the above two types of applications of Azuma's inequality considering both the optimization problem and the random graph factors. Compared with W.T.Rhee and M. Talagrand's work \cite{Rhee1987}, the stochastic model is replaced by a random network, and instead of an optimal solution sequence of increasing system sizes, we consider an arbitrary solution on a fixed but gradually uncovered random graph and finally try to find the optimal ones. Shamir, Spencer's \cite{Shamir1987} and C. McDiarmid's \cite{McDiarmid1989} works are the most relevant to ours. We use almost the same technique in derivation. Other works mentioned above all consider some direct statistics of random networks without an optimization procedure so is conceptually simpler than ours.

\section{General Model}
The random network model is represented by an ensemble $\Omega(N)=\{\Gamma(N),\mathscr{P}\}$, where $\Gamma(N)$ is the collection of all connected graphs with N nodes and $\mathscr{P}$ is a probability measure on $\Gamma(N)$. $G\in \Omega(N)$ is considered as a random network taking
values in $\Gamma(N)$ with probability $\mathscr{P}(G)$. Its adjacency matrix is given by $A=[A_{ij}]_{N\times N}$ and the degree of each node is denoted by $d_i$ ($i=1,...,N$). Additionally we require the probabilities for any two nodes to be linked are independent, i.e.
$\{A_{ij}|i<j\}$ are $N(N-1)/2$ independent random variables. A spin vector $\vec{s}=(s_1,...,s_N)$ $(s_i=\{-1,1\})$ assigns a spin to each node and indicates a partition which takes the nodes with the same spin in one group. The community detecting function is given by:

\begin{equation}
h_G(\vec{s})=-{1\over N}\sum_{i,j}F_{ij}(G)s_i s_j
\end{equation}

where $\{F_{ij}\}$ are functions of the random network $G$ and hence random variables. Because the form of community detecting function can be considered as the Hamiltonian of a spin system, we also call $h_G(\vec{s})$ the Hamiltonian in later context. Let $S$ be the spin configuration space including all possible spin states satisfying given constrains. We take $\vec{s}_0\in S$ as the optimal configuration of $h_G$, and

$$H(G)=\min_{\vec{s}\in S}{h_G(\vec{s})}=h_G(\vec{s}_0)$$

In stead of $h_G(\vec{s})$, we focus on $H(G)$ because this value only depends on $G$ therefore is a property of the network itself rather than some specific configuration. Furthermore the statistics of $H(G)$ provides information of the network ensemble. In this section, we give large deviation results of the $H(G)$ distribution.

%In the following section, we will prove concentration results for
%several community structure (C.S.) and graph partition problems, for
%specifically, the non-sparse and sparse bounded degree subclasses of
%ER and Chung-Lu. For all the non-sparse cases, we will use Theorem 1
%together with a ``divide and conquer" technique to obtain
%concentration. For the sparse bounded degree subclasses, we will use
%Theorem 2 along with technical estimates specific to each problem.

{\bf Theorem 1}  \emph{If $\{F_{ij}\}$ ($i,j=1,...,N$) satisfy that for any two networks $G$,$G'\in\Omega$ which only differ by one edge, $|F_{ij}(G)-F_{ij}(G')|\leq c$
where $c$ is a constant independent of $i,j$, then for every real number $t$, }

\begin{equation}
\mathscr{P}(|H(G)-<H>_\Omega|>t)\leq 2 \exp{(-{t^2\over c^2})}
\end{equation}

\emph{where $<\cdot>_\Omega$ is the ensemble average and $c^2$ is a constant independent of $t$ and $\Omega(N)$.}\\

To prove theorem 1, we need Azuma's inequality and the following lemma:

{\bf Theorem 2(Azuma)} \cite{Azuma} \emph{Suppose $Y_K=\sum_{i=1}^K X_i$ is a martingale, given boundedness of each increment $|X_i|\leq b_i$, the inequality holds for any
real number $t$:}
$$\mathscr{P}(|Y_K-E[Y_K]|>t)\leq 2\exp({-{t^2\over 2\sum_i b_i^2}}).$$\\

%%%McDiarmid's Inequality

{\bf Lemma 1} \emph{If the conditions in Theorem 1 hold, then $|H(G)-H(G')|\leq {c\over N}$.}

{\bf Proof of Lemma 1:}
At first, suppose $G$ and $G'$ only differ at the edge between the nodes $i_0$ and $j_0$.
$$|h_G(\vec{s})-h_{G'}(\vec{s})|\leq {1\over N} |F_{i_0j_0}(G)-F_{i_0j_0}(G')||s_0 s_0|\leq {c\over N}$$
Without loss of generality, we assume $H(G)\geq H(G')$ and $H(G)$ is achieved at a specific configuration $\vec{s_0}$, that is,
$h_G(\vec{s_0})=H(G)\geq H(G')\geq h_{G'}(\vec{s_0})$. So $|H(G)-H(G')|\leq|h_G(\vec{s_0})-h_{G'}(\vec{s_0})|\leq{c\over N}$.\\

{\bf Proof of Theorem 1}
Suppose $G\in \Omega(N)$ is a random graph whose node pairs are labeled from $1$ to $N(N-1)/2$. Each node pair corresponds to one random variable $A_k=\bold{1}_{\{\text{node pair $k$ is linked}\}}$ $(k=1,...,{N(N-1)/2})$. Therefore the random graph $G$ can be considered as a random process $\{A_1,...,A_{N(N-1)/2}\}$, and its filtration is $F_k$ $(k=0,...,{N(N-1)/2})$.
Define $H_k=E[H(G)|F_k]$ $(k=1,...,{N(N-1)/2})$, and $H_0=E[H(G)|F_0]=<H>_\Omega$. $H_k$ $(k=0,...,{N(N-1)/2})$ is a martingale by construction. We construct an auxiliary process $G'=\{A'_1,...,A'_{N(N-1)/2}\}$ for which $A'_j=A_j$ when $j\neq k+1$ and $A'_j=1-A_j$ when $j=k+1$. $G'$ shares the same filtration with $G$ and represents the graph that only differs from $G$ by the link between the node pair $k+1$.
The increment of $H_k$
\begin{eqnarray*}
|H_{k+1}-H_k|&=& \left| H_{k+1}-E[H_{k+1}|F_k] \right| \\
             &=& \left| E[H(G)|F_{k+1}]-E[H(G')|F_{k+1}]\right|\\
             &\leq& E[|H(G)-H(G')| | F_{k+1}]\\
             &\leq& {c\over N}
\end{eqnarray*}
In the last inequality, $|H(G)-H(G')|$ as a random variable only depending on the realization of $G$ is bounded by ${c\over N}$ according to Lemma 1, so its conditional expectation is also controlled by this bound. We apply Azuma's inequality on $H_k$, and obtain:
$$\mathscr{P}(|H(G)-<H>_\Omega|>t)\leq 2 \exp{(-{t^2\over 2 ({c\over N})^2 N(N-1)/2})}\leq 2 \exp{(-{t^2\over c^2})} \quad \bigbox$$

The inequality by itself only indicates some level of stability of $H(G)$ over the ensemble $\Omega$. To interpret this inequality as a concentration result, instead of a single network ensemble, we need a random network model consisting of a sequence of network ensembles with different system size $N$.

The above general model does not specify the way to generate the ensemble sequence when $N\go \infty$, especially the growing speed of the number of edges. To keep the networks connected, the number of edges must grow at least linearly with respect to the number nodes. We take the network models whose numbers of edges grows linearly and super-linearly as the sparse networks and non-sparse networks separately. There is no general concentration result for the non-sparse networks, so we leave the discussion of this case in the later sections for specific community detecting functions. As to the sparse networks, there are various ways to define them as long as their numbers of edges grow linearly and we only consider one definition, K-bound networks.

\section{K-bound networks}
A given network is called K-bound if the constant $K$ is an upper bound of the network degrees. A network ensemble or a random network model is called K-bound if $K$ is a uniform upper bound of the degrees over the network ensemble or the ensemble sequence.

A subtle point here is that the K-bound constrain may violate the independence of the links. But since in a random network with independent links, the probability distribution of the degree decays very fast (the probability for the degree exceed $K$ has the order $O(e^{-K^2})$), when $K$ increase, the difference between the probability measures of the network ensembles with and without the K-bound constrain goes to zero very quickly. For K-bound networks, we have the concentration result:

{\bf Theorem 3} \emph{If $\Omega$ is a K-bound ensemble, $\{F_{ij}\}$ ($i,j=1,...,N$) satisfy that for any two networks $G$,$G'\in\Omega$ which only differ by one edge, $|F_{ij}(G)-F_{ij}(G')|\leq c$ where $c$ is a constant independent of $i,j$, then for every real number $t$,}
\begin{equation}
\mathscr{P}(|H(G)-<H>_\Omega|>t)\leq 2 \exp{(-{Nt^2\over 8K^2c^2})}
\end{equation}

{\bf Proof of Theorem 3} Consider $G\in \Omega(N)$ is a random graph labeling the nodes as $1,2,...,N$. $G_k$ is the subgraph containing nodes from $1$ to $k$ and also considered as a filtration.
Define $H_k=E[H(G)|G_k]$ $(k=1,...,N)$, and $H_0=E[H(G)|G_0]=<H>_\Omega$. $H_k$ $(k=0,...,N)$ is a martingale by construction. Consider $G''$ is a network only differ from $G$ by one node $k+1$, ie. all the different elements in the adjacency matrix are for the edges linked to node $k+1$. $G''_k$ is the corresponding filtration for $G$. Note that the degree is at most $K$, therefore $G$ and $G''$ can only differ by at most $2K$ edges.
\begin{eqnarray*}
|H_{k+1}-H_k|&=& \left| H_{k+1}-E[H_{k+1}|F_k] \right| \\
             &=& \max\{\left| E[H(G)|G_k]-E[H(G'')|G''_k]\right|\}\\
             &\leq& 2K\left| E[H(G)|G_k]-E[H(G')|G_k]\right|\\
             &\leq& {2Kc\over N}
\end{eqnarray*}

Here $G'$ is the same as defined in the proof of Theorem 1. The difference $H(G)-H(G'')$ can be decomposed into at most $2K$ terms which all have the same bound as $H(G)-H(G')$. Apply Azuma's inequality, we have:
$$\mathscr{P}(|H(G)-<H>_\Omega|>t)\leq 2 \exp{(-{t^2\over 2 ({2Kc\over N})^2 N})}\leq 2 \exp{(-{Nt^2\over 8K^2c^2})} \quad \bigbox$$

{\bf Remark:} According to this inequality, the fluctuation of $H(G)$ about its mean has the order $O(N^{-\half})$.

When we are given a network $G_0=\{A_{0ij}\}$, one way to generate a network ensemble is to consider the given network undergoes a certain perturbation. Suppose the perturbation is denoted by $\delta G(p_0,p_1)$ where $p_0,p_1$ are two probabilities. The network ensemble $G=\{A_{ij}\}$ generated by $G_0$ and $\delta G(p_0,p_1)$ satisfies:\\
(a)if $A_{0ij}=1$, then $A_{ij}=0$ with probability $p_0$ and $A_{ij}=1$ with probability $1-p_0$.\\
(b)if $A_{0ij}=0$, then $A_{ij}=1$ with probability $p_1$ and $A_{ij}=0$ with probability $1-p_1$.

To conserve the average degree of the original network $G_0$, we additionally require $p_0 m=p_1 ({N(N-1)\over 2}-m)$, where $m$ is the total number of edges in $G_0$. Since $m<{NK\over 2}$, the requirement implies:
$$p_1\leq {K\over N-1-K}p_0$$

If $G_0$ satisfies the K-bounded degree condition, we have:

{\bf Theorem 4} \emph{If $\Omega$ is a K-bound ensemble generated by a given K-bound network $G_0$ and a small perturbation $\delta G$, $\{F_{ij}\}$ ($i,j=1,...,N$) satisfy that for any two networks $G$,$G'\in\Omega$ which only differ by one edge, $|F_{ij}(G)-F_{ij}(G')|\leq c$ where $c$ is a constant independent of $i,j$, then for every real number $t$,}
\begin{equation}
\mathscr{P}(|H(G)-<H>_\Omega|>t)\leq 2 \exp{(-{Nt^2\over 4[c^2(2K^2p_0^2+Kp_0)+Kct/3]})}
\end{equation}

To prove Theorem 4, we need a variance form of Azuma's inequality as below:

{\bf Theorem 2'} \cite{Chung2006} \emph{Suppose $Y_K=\sum_{i=1}^K X_i$ is a martingale, given boundedness of each increment $|X_i|\leq M$, $Var(X_i)= v_i^2$ the inequality holds for any
real number $t$:}
$$\mathscr{P}(|Y_K-E[Y_K]|>t)\leq 2\exp({-{t^2\over 2(\sum_i v_i^2+Mt/3)}})$$.\\

{\bf Proof of Theorem 4}
The proof is almost the same as that for Theorem 3. We only need to replace Theorem 2 in that proof by Theorem 2' and find the upper bound of $Var(X_i)$.

Let random variable $L_i$ be the number of changes between $G_0=\{A_{0ij}\}$ and $G=\{A_{ij}\}$ $(i<j)$ related to the node $i$. Since $L_i$ can be considered as a sum of 0-1 random variables:
\begin{eqnarray*}
Var(L_i)&=&K p_0(1-p_0)+ (N-1-K)p_1(1-p_1)\\
        &\leq& K p_0+ (N-1-K)p_1\\
        &\leq& 2Kp_0\\
\end{eqnarray*}
and
$$E[L_i]=Kp_0+(N-1-K)p_1\leq 2Kp_0$$

Since $|X_i|\leq L_i{c\over N}$, we have:
$$Var(X_i)\leq{c^2\over N^2} E[L_i^2]={c^2\over N^2}(E^2[L_i]+Var(L_i))={c^2\over N^2} (4K^2p_0^2+2Kp_0),$$
and Theorem 2' becomes the conclusion of Theorem 4. Compared with Theorem 3, Theorem 4 gives a much sharper concentration around the average value whenever $p_0$ is small. $\quad \bigbox$

%According to Theorem 4, we conclude that even if a network ensemble is not exact K-bound, whenever it can be considered as a K-bound network under a small perturbation, the value of community detecting function on this ensemble is still stable and concentrated.

In the following sections, we will prove the concentration for several different community detecting functions on both the non-sparse and sparse subclasses of Erdos-Renyi (ER) and Chung-Lu (CL) random networks. For all the non-sparse cases, we can use Theorem 1 together with a ``divide and conquer'' technique described in the next section to obtain concentration. However, when the community detecting function have additional properties like $<H(G)>_\Omega \simeq O(N)$ as in the graph bipartitioning case, we can scale $t$ in Theorem 1 by $N$ to obtain a better concentration inequality without using the ''divide and conquer'' technique. For the K-bounded degree subclasses of ER and CL networks, we will use Theorem 3/Theorem 4 along with technical estimates specific to each case.

\section{Modularity}
%Modularity is one of the criterion for community detection \cite{7}.
%In this case $F_{ij}(G)={N\over 4m} (A_{ij}-{d_i d_j\over 2m})$, therefore
%the Hamiltonian (1) is
Modularity is one of several effective criterion for the detection
of community structures in random networks \cite{7}. In this
section, the above theorems are used to obtain LDP results on the
Modularity functional over the non-sparse Erdos-Renyi subclass
$ER[p,N]$ and where there is a uniform bound on the degrees of $G,$
over the sparse $ER[Np,K]$ subclass (and Chung-Lu $CL[N,\varpi ,K]$
subclass which properties are given in the next section). The
non-sparse cases here and in the next section are treated by the
''Divide and Conquer'' technique because, without prior knowledge on
properties of their respective ensemble average, $<H(G)>_\Omega ,$
we have no recourse to the shorter method of rescaling $t$ in
theorem 1 by $N.$

In the Modularity problem, $F_{ij}(G)={\frac N{4m}}(A_{ij}-{\frac{d_id_j}{2m}%
})$, therefore the Hamiltonian (1) is
\begin{equation}
h_G(\vec{s})=-{1 \over 4m}\sum_{i\neq j} (A_{ij}-{d_i d_j\over 2m}) s_i s_j
\end{equation}
Here $m$ is the total number of edges in the network, $d_i$ is the degree of node $i$, and $A_{ij}$ is the adjacency matrix.
The Hamiltonian is more complicated than the previous two cases as $F_{ij}(G)$  depends on not only the local information $A_{ij}$ but also the global information
such as $d_i, m$. Suppose $G'$ only have one more edge ($i_0j_0$) than $G$.
To apply Theorems 1 and 3, , we need to estimate $|F_{ij}(G)-F_{ij}(G')|$ as follows,
\begin{eqnarray}
&&{4m\over N}|F_{ij}(G)-F_{ij}(G')| \notag \\
&=&\left|\sum_{i\neq j} \left[(A'_{ij}-A_{ij})-({d'_i d'_j\over 2(m+1)}-{d_i d_j\over 2m})\right] s_i s_j\right| +({4m\over N})({m+1\over m}-1)|F_{ij}(G')|\notag\\
&\leq&{1\over m} |A_{i_0,j_0}-A'_{i_0,j_0}|+\left|\sum_{i\neq j}{d'_i d'_j-d_i d_j\over 2(m+1)}\right|+\left|\left({1\over 2(m+1)}-{1\over 2m}\right)\sum_{i\neq j}{d_i d_j}\right|\notag\\
&+&({1\over m+1})\left|A'_{ij}-{d_id_j\over 2(m+1)}\right|
\end{eqnarray}\\

Next, we estimate the 4 terms in equation (6). The first term $|A_{i_0,j_0}-A'_{i_0,j_0}|=1$. The second term:
\begin{eqnarray*}
&&\left|\sum_{i\neq j}{d'_i d'_j-d_i d_j\over 2(m+1)}\right|\\
&=&{1\over 2(m+1)} \left|\sum_{j\neq i_0}(d'_{i_0}-d_{i_0})d_j + \sum_{i\neq j_0}(d'_{j_0}-d_{j_0})d_i + (d'_{i_0}-d_{i_0})(d'_{j_0}-d_{j_0})\right|\\
&\leq& {1\over 2(m+1)}(2m+ 2m +1)\\
&<& 2
\end{eqnarray*}

For the third term:
\begin{eqnarray*}
&&\left|\left({1\over 2(m+1)}-{1\over 2m}\right)\sum_{i\neq j}{d_i d_j}\right|\\
&\leq&{1\over 2m(m+1)}|(\sum_i d_i)(\sum_j d_j)|\\
&\leq&{1\over 2m(m+1)}(2m)^2\\
&<& 2
\end{eqnarray*}

For the last term, since $d_i,d_j\leq m$,
\begin{eqnarray*}
&&({1\over m+1})\left|A'_{ij}-{d_id_j\over 2(m+1)}\right|\\
&\leq&({1\over m+1})\left|{m^2\over 2(m+1)}-1\right|\\
&\leq& \max\{ {m^2\over 2(m+1)^2},{1\over m+1}\}\\
&<& {1\over 2}
\end{eqnarray*}

With the above inequalities, we conclude
\begin{equation}
|F_{ij}(G)-F_{ij}(G')|\leq {11N \over 8m}\leq {11N\over 8m^*}
\end{equation}
where $m^*$ is the minimum value of $m$ required by connectivity. This completes the proof of the technical estimates needed in the application of Theorems 1 and 3 below.

%For general network without additional assumptions, according to Theorem 1, we have:
For general ER networks without additional non-sparseness or bounded
degree properties and also for sparse ER networks without
assumptions on for instance the scale of the ensemble average,
$<H(G)>$, Theorem 1 gives:
$$\mathscr{P}(|H(G)-<H>_\Omega|>t)\leq 2 \exp{(-{t^2\over \sigma^2 ({N\over m^*})^2 })}\leq2 \exp{(-{t^2\over \sigma^2})}$$
where $\sigma={11\over 8}$, which shows that the direct application of Theorem 1 alone does not allow us to obtain concentration.

However for sparse ER network with bounded degree, using Theorem 2, we prove:

{\bf Theorem 5} \emph{In $ER[Np=\lambda]$ and $CL[N,\vec{w}]$ network ensemble $\Omega$ with uniform upper bound $K$ of the degree, the optimal modularity $H(G)$ satisfies a concentration inequality:}
\begin{equation}
\mathscr{P}(|H(G)-<H>_\Omega|>t)\leq 2 \exp{(-{Nt^2\over {25\over 2}({N\over m^*})^2 K^2})}\leq 2 \exp{(-{Nt^2\over {25\over 2}K^2})}
\end{equation}

In some sense, the networks with degree upper bound is an extreme
case of sparse network. Next, we consider the other extreme case,
the non-sparse networks, using the following ``divide and conquer" method. The network is non-sparse if:\\
(a). $<m>_{\Omega(N)}=pN^\alpha/2$ ($p,\alpha$ are two constant parameter)\\
(b). $P(m-<m>_{\Omega(N)}<-N t)\leq \exp{(-{t^2\over \lambda^2})}$\\
where $m$ is the total number of edges, $p$ is a constant independent of $N$.
The meaning of property (b) is shown in the following lemma:

{\bf Lemma 2} \emph{If in a network all the edges are independent, there exist a constant $\lambda$ independent of $N$, s.t. the number of edges $m$ satisfies}
$$P(m-<m>_{\Omega(N)}<-N t)\leq \exp{(-{t^2\over \lambda^2})}$$
{\bf Proof of Lemma 2}:
$m=\sum_{i<j} \bold{1}_{\{A_{ij}=1\}}$, where $\bold{1}_{\{A_ij=1\}}\in[0,1]$ are independent random variables.
According to Hoeffding inequality \cite{Hoeffding1963},
$$\mathscr{P}({m-<m>\over N(N-1)/2}<-t)\leq \exp(-{2t^2 N^2(N-1)^2\over 4\sum_{i<j}1})$$
Substituting $t$ by $2t/(N-1)$, we get
$$\mathscr{P}(m-<m><-Nt)\leq \exp(-{2t^2 N^2\over N(N-1)/2})\leq \exp(-{t^2 \over \lambda^2})$$
where $\lambda={1\over 2}$.

An example of non-sparse network is ER network $ER(p,N)$ with constant probability $p$ for any two node to be linked.
Replacing $t$ in (b) by $\mu N^{\alpha-1}/2$, where $\mu$ is independent of $N$, we get
\begin{eqnarray}
\mathscr{P}(m< (p-\mu) N^\alpha/2) &\leq& \mathscr{P}(m< (<m> _{\Omega(N)}-\mu N^\alpha/2))\nonumber\\
                     &\leq& \exp{(-{\mu^2 N^{2(\alpha-1)}})}
\end{eqnarray}

Finally, we split $P(|H(G)-<H>_\Omega|>t)$ into two cases according to $m$, and get the concentration result:
\begin{eqnarray}
&&P(|H(G)-<H>_\Omega|>t) \nonumber \\
&=&P(m< (p-\mu) N^\alpha/2)P( |H(G)-<H>_\Omega|>t |m< (p-\mu) N^\alpha/2) \nonumber \\
&+& P(m\geq (p-\mu) N^\alpha/2)P(|H(G)-<H>_\Omega|>t |m\geq (p-\mu) N^\alpha/2)\\
&\leq& \exp{(-{\mu^2 N^{2(\alpha-1)}})} \cdot 2 \exp{(-{t^2})} +1 \cdot 2 \exp{(-{t^2\over ({N\over (p-\mu)N\alpha})^2 })} \nonumber \\
&\leq& 2 \exp{(-{\mu^2 N^{2(\alpha-1)}}-{t^2})}+2 \exp{(-{t^2 N^2(p-\mu)^2})}
\end{eqnarray}

When $1<\alpha\leq 2$, the inequality shows concentration. For ER[N,p] network, $\alpha=2$, we have the following theorem.

{\bf Theorem 6} \emph{In non-sparse ER[N,p] network ensemble $\Omega$, the optimal modularity $H(G)$ satisfies the concentration inequality:}
\begin{equation}
P(|H(G)-<H>_\Omega|>t)\leq 2 \exp{(-{\mu^2 N^2}-{t^2})}+2 \exp{(-{t^2 N^2(p-\mu)^2})}
\end{equation}

\section{ Concentration of Modularity on Chung-Lu network $CL_\infty [N,\varpi
,\beta ]$, $CL[N,\varpi ,K]$}

Scale free or power law random graphs including Barabasi-Albert
(BA), Molloy-Reed (MR) and the Chung - Lu (CL) models are frequently
encountered in the study of random networks that arise in social and
ecological problems. For brevity, we will use the Chung-Lu model
which is based on
fixed expected degrees sequence $\varpi $ connected random graphs $%
CL(N,\varpi )$ \cite{CL}. The Chung-Lu model is easier to work with than
the Molloy-Reed \cite{MR}, Newman-Strogatz-Watts \cite{WS} and Barabasi-Albert \cite{BA}
formulations because it specifies more information at the level of
each node $j$ in the graph $G(N)$. It is based on working with
subsets $CL[N,\varpi ]$ of the scale free random graphs specified by
a fixed (deterministic) expected degrees sequence of weights
\begin{eqnarray*}
\varpi &=&(w_1,...,w_N) \\
w_j &=&E[\deg (j)].
\end{eqnarray*}
The average degree in $CL[N,\varpi ]$ is given by
\[
d(\varpi )=\frac 1N\sum_{l=1}^Nw_l.
\]

In the growth process of CL random networks (cf. BA \cite{BA}), a new node
$v_i$ is added at time $i\leq N$ and $m^{\prime }$ random and
independent edges are then added between this node $v_i$ and those
already present. Thus, for node $i$ added, the probability of adding
a link to node $j$ is $w_j/\sum_lw_l.$ The second order degree or
average number in the second generation of nodes (where the random
graph connectedness problem is viewed as a two stage branching
process ) is given by
\[
\bar{d}(w)=\sum_{j=1}^nw_j\frac{w_j}{\sum_lw_l},
\]
whence it is easily shown via an application of Cauchy-Schwartz that
\[
\bar{d}(\varpi )-d(\varpi )\geq 0.
\]

A key property of $CL[N,\varpi ]$ for the proof of concentration of
Modularity below is the probability for an edge between arbitrary
pair of nodes $i$ and $j$ or independent random variables $A_{ij}=1$
(where $A(G)$ is the adjacency matrix of the random graph $G)$ is
given by \cite{Chung2002a,Chung2002b}
\[
\Pr \{A_{ij}(G(N))=1\}=\frac{w_iw_j}{\sum_{l=1}^Nw_l}.
\]
In each subclass $\Omega =CL[N,\varpi ],$ the average number of
edges $m$ is

\begin{eqnarray*}
&<&m>_\Omega =\left( \frac{N(N-1)}2+N\right) \frac 1{N^2}\sum_i^N\sum_j^N%
\frac{w_iw_j}{\sum_kw_k} \\
&=&\frac{N(N+1)}{2N^2}\sum_i^Nw_i=\frac{d(\varpi )}2(N+1),
\end{eqnarray*}
where the average degree $d$ may depend on $N$ through the weights
$w_j.$

There are no known concentration results on the whole class
$CL[N,\varpi ]$; this class includes sparse random networks that are
not uniformly bounded in degree of nodes. Using the ''divide and
conquer'' technique and theorem 1 twice, we prove LDP results for
non-sparse Chung-Lu networks $CL_\infty [N,\varpi ,\beta ]$ for
which in addition to the above properties, the average degree grows
with $N,$ that is, for $\beta >0$ and $B>0,$ both
independent of $N,$%
\[
d(\varpi )\geq BN^\beta .
\]
Concentration of Modularity in the subclass of bounded degree
Chung-Lu networks $CL[N,\varpi ,K]$ is proved in the previous
section.

Next by Lemma 2, this subclass $\Omega (N)=CL_\infty [N,\beta ]$
satisfies
the non-sparse property (b) in :\\(a). $<m>_{\Omega (N)}\geq MN^\alpha $, ($%
M>0$, $\alpha >1$ are two constant parameters).\\(b) there exists constant $%
\lambda $ independent of $N$ such that $\Pr (m-<m>_{\Omega
(N)}<-Nt)\leq \exp {(-{\frac{t^2}{2\lambda ^2}})}$ for any real
$t>0$ $.$

Property (a) holds by choosing $M=B/2,$and $\alpha =1+\beta $.
Replacing $t$ in (b) by $\mu N^\beta $ gives
\begin{eqnarray*}
P(m<(M-\mu )N^{\beta +1}) &\leq &P(m<(<m>_\Omega -\mu N^{\beta +1})) \\
&\leq &\exp {(-{\frac{\mu ^2N^{2\beta }}{2\lambda ^2}}).}
\end{eqnarray*}
By conditioning $P(|H(G)-<H>_\Omega |>t)$ according to $m$, we
derive the result: for constants $\lambda ,\beta >0$ and $\mu
<M=B/2,$ that do not depend on $N,$
\begin{eqnarray*}
&&P(|H(G)-<H>_{CL_\infty }|>t) \\
&=&P(m<(M-\mu )N^{\beta +1})P\left\{ |H(G)-\text{ }<H>_{CL_\infty
}|>t\text{
}|\text{ }m<(M-\mu )N^{\beta +1})\right\} \\
&+&P(m\geq (M-\mu )N^{\beta +1})P\left\{ |H(G)-<H>_{CL_\infty }|>t\text{ }|%
\text{ }m\geq (M-\mu )N^{\beta +1})\right\} \\
&\leq &\exp {(-{\frac{\mu ^2N^{2\beta }}{2\lambda ^2}})}\cdot 2\exp {(-{%
\frac{t^2}{2\sigma ^2}})}+1\cdot {2\exp }\left( -\frac{t^2(M-\mu
)^2N^{2\beta }}2\right) \\
&\leq &2\exp {(-{\frac{\mu ^2n^{2\beta }}{2\lambda
^2}}-{\frac{t^2}{2\sigma ^2}})}+{2\exp }\left( -\frac{t^2(M-\mu
)^2n^{2\beta }}2\right) .
\end{eqnarray*}
We used the following applications of theorem 1 : for any real $t>0,$%
\begin{eqnarray*}
P\{|H(G)-\text{ } &<&H>_{CL_\infty }|>t\text{ }|\text{ }m<(M-\mu
)N^{\beta
+1}\}\leq 2\exp {(-{\frac{t^2}{2\sigma ^2}})} \\
P(|H(G)- &<&H>_{CL_\infty }|>t\text{ }|\text{ }m\geq (M-\mu
)N^{\beta +1}\}\leq {2\exp }\left( -\frac{t^2(M-\mu )^2N^{2\beta
}}2\right)
\end{eqnarray*}
with respectively,
\begin{eqnarray*}
c^2\frac{N(N+1)}2 &=&{\frac{25N(N+1)}{32m^2}} \\
&\leq &{\frac{25}{32}}(1+{\frac 3{N-1}+}\frac 2{(N-1)^2}) \\
&\leq &{\frac{25}{16}\equiv \sigma }^2,
\end{eqnarray*}
and when $m^2\geq (M-\mu )^2N^{2(\beta +1)},$%
\begin{eqnarray*}
c^2\frac{N(N+1)}2 &=&{\frac{25N(N+1)}{32m^2}} \\
&\leq &\frac{25}{32}\frac{N(N+1)}{(M-\mu )^2N^{2\beta +2}} \\
&\leq &\frac 1{(M-\mu )^2N^{2\beta }},
\end{eqnarray*}
where the constant $c>0$ comes from the technical estimate (6) that
is valid for any pair of random networks $G,G^{\prime }$ differing
by exactly one edge, in any random ensemble $\Omega .$ This
completes the proof of

\textbf{Theorem 7} \emph{In non-sparse }$CL_\infty $\emph{[}$N,B,\beta $%
\emph{] network ensemble $\Omega $, the optimal modularity $H(G)$
satisfies the concentration inequality: there exists} $\mu >0$
independent of $N,$ such that for any real $t>0,$
\begin{equation}
P(|H(G)-<H>_\Omega |>t)\leq 2\exp {(-{\frac{\mu ^2N^{2\beta }}{2\lambda ^2}}-%
{\frac{t^2}{2\sigma ^2}})}+{2\exp }\left( -\frac{t^2(B/2-\mu )^2N^{2\beta }}%
2\right)
\end{equation}
where $\lambda =1/2,$ and $\sigma =5/4.$

\section{Graph Bipartitioning on ER and CL Networks}
The so-called graph bipartitioning problem \cite{5} is the simplest example where $F_{ij}(G)=A_{ij}$ in (1).
\begin{equation}
h_G(\vec{s})=-{1\over N}\sum_{i,j}A_{ij}s_i s_j
\end{equation}
In this problem, each spin state corresponds to a two group partition of the given graph, and the optimum gives the partition with the least intergroup links.
Fu and Anderson \cite{Fu1986} has shown the equivalence between this problem for $ER[N,p]$
network (any two nodes are linked by an edge with probability $p>0$ where $p$ is independent of $N$)
and the infinite range SK model and derived the average solution in the thermodynamic limit.
Banavar et al. \cite{Banavar1987} investigate another case of $ER[Np=\lambda]$ network and give
an empirical average solution. However, without a concentration result, these solutions are
only heuristic. Even if we accept the solution in the thermodynamic limit, the errors of the
solutions applied to a finite system cannot be estimated. Our research complete their results.

In their studies, the constrain of zero magnetization $M=\sum_i
s_i=0$ is required which force an equal size partition. With a
non-zero fixed magnetization constrain $\sum_i s_i=c\neq 0$, we get
the optimal partition with given two group sizes. Without any
constrain on $M$, we get the overall optimal partition considering
all possible group sizes. No matter which constrain we use, it only
changes the spin configuration space and hence the definition of
$H(G)=\min_{s\in S}h_G(s)$, so we have exactly the same
concentration result and proofs in all of these cases.

For the most general random graph ensemble, we apply Theorem 1 as follows. Since
\begin{equation}
|F_{ij}(G)-F_{ij}(G')|=|A_{ij}(G)-A_{ij}(G')|\leq 1
\end{equation}
we have
$$\mathscr{P}(|H(G)-<H>_\Omega|>t)\leq 2 \exp{(-t^2)}$$
which shows that concentration is not obtained by this approach. To
compare with the papers mentioned above, we consider the two types
of ER network separately. For the first type $ER[N,p]$, we assume in
addition, $<H>_\Omega$ is of the order $O(N)$ which is consistent
with the results from the replica method \cite{Fu1986}, which
allows us to scale $t$ with the same order and have:

{\bf Theorem 8} \emph{In a non-sparse $ER[N,p]$ network ensemble
$\Omega$, the optimal graph partitioning $H(G)$ satisfies the
concentration inequality:}
$$P(|H(G)-<H>_\Omega|>Nt)\leq 2 \exp{(-N^2t^2)}$$

For the second type $ER[Np=\lambda]$, for which $Np$ is fixed as $N\go \infty$, we prove a theorem under the additional assumption, bounded degree. Consider the $ER[Np=\lambda]$ network ensemble generated with parameter $N,p$ excluding all samples whose maximum degree exceed the bound $K>Np$, where $K$ is independent of $N$. So when $N\go \infty$, the degree distribution of the network tends to a Poisson distribution with expected degree $Np$ but with a cutoff at $K$. If $K$ is large enough, the network ensemble generated like this is almost the second type of ER. For this case, we prove use theorem 3 to prove:

{\bf Theorem 9} \emph{In $ER[Np=\lambda]$ ensemble $\Omega$ with
uniform upper bound $K$ of the degree, the optimal graph
partitioning $H(G)$ satisfies the concentration inequality:}
$$\mathscr{P}(|H(G)-<H>_\Omega|>t)\leq 2 \exp{(-{Nt^2\over 8 K^2})}$$

For comparison with Theorems 8 and 9, the simulated concentration for both cases of ER model
is shown in the following figures. For the constant $p$ case, the fluctuation
of $H(G)/N$ is roughly of the order $O(N^{-1})$. While for the constant $Np$
case, the fluctuation of $H(G)$ is roughly of the order $O(N^{-1/2})$.

\begin{center}
\begin{figure}[!htbp]
  \includegraphics[width=0.8\textwidth]{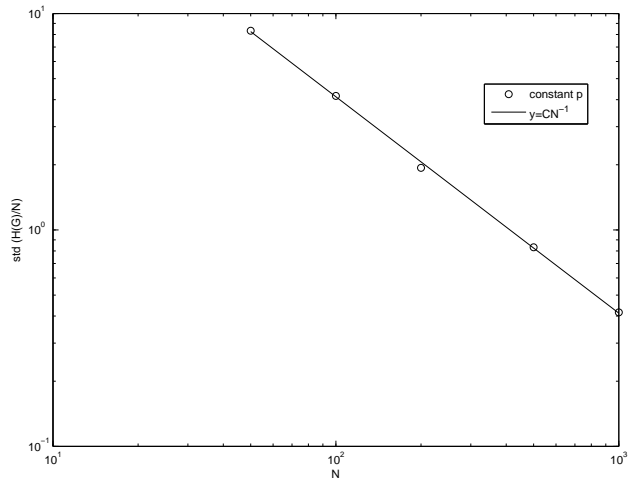}\\
\includegraphics[width=0.8\textwidth]{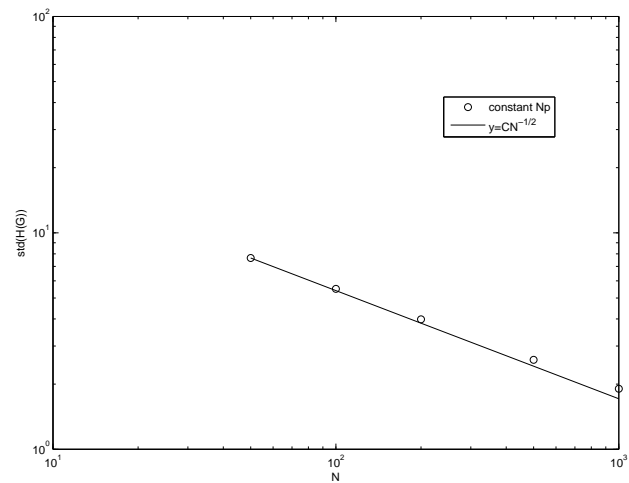}\\
  \caption{Each data point is obtained by optimization using simulated annealing on 100 realization of ER network. The first figure are for ER with constant $p=0.05$, and the second figure are for ER with constant $pN=5$}
  \label{figure:1}
\end{figure}
\end{center}

%\section{relationship to SK model, Sample average and quenched average, phase transition}
As an application of the above theorems to a real world problem, we consider
the circuit partition optimization problem proposed by S. Kirkpatrick et al \cite{8}.
The objective function they proposed is equivalent to:
$$h_G(\vec{s})= -{1\over N}\sum_{i,j}(A_{ij}/2-\lambda)s_i s_j$$
So $F_{ij}=A_{ij}-\lambda$, where $A_{ij}$ is the number of signals passing between circuits $i$ and $j$, $\lambda$
is a balancing coefficient in optimization \cite{8}. %\lambda<z/2
In their study, $A_{ij}$ is no longer a $0-1$ random variable, but
as long as the $A_{ij}$'s are uniformly bounded which is quite
reasonable in their problem, we can always normalize $A_{ij}$ by the
uniform upper bound and all the theorems and proofs in this section
still hold. The $\lambda$ introduced in this problem represent a
penalty on partitions with nonequal groups sizes, and is a special
case of that in the modularity function.

\section{Concentration for Q-Potts Community Structures on ER and CL networks
}

This section is devoted to concentration results for the family of
objective functionals derived from q-Potts models which are
introduced as a viable alternative to the Modularity community
detection algorithm:
\begin{eqnarray}
h_G(\vec{s})=-\frac J{4m}\sum_{i\neq j}A_{ij}\delta (s_i,s_j)+\frac
\gamma {2m}f(\vec{n}_s) \label{h}
\end{eqnarray}

where the spins $s_j\in (0,1,...,q-1),$ $A_{ij}=0,1$ is the
adjacency matrix
of $G,$ $\delta (s_i,s_j)=1$ only if $s_i=s_j,$ otherwise zero, $f(\vec{n}%
_s) $ is a function of the occupation numbers $\vec{n}%
_s=(n_0,n_1,....,n_{q-1}),$ $G(N,m)$ is a random graph with $N$
nodes and $m$ links, $J$ is a ferromagnetic interaction energy, and
$\gamma $ is a parameter that determines the antiferromagnetic
activity of the Hamiltonian$. $ The particular example of these
Hamiltonians studied by Reichardt\cite{Reichardt2004} is
\[
h_G(\vec{s})=-\frac J{4m}\sum_{i\neq j}^NA_{ij}\delta
(s_i,s_j)+\frac \gamma {2m}\sum_{s=0}^{q-1}\frac{n_s(n_s-1)}2
\]
where for $\gamma >\gamma ^{*},$ the optimum of $h_G$ favors
community structures that reflect the network topology of $G.$
Threshold $\gamma ^{*}$ is fixed by requiring $h_G(homogeneous)>h_G(diverse)$ which for two communities $c(n_1,m_1)$ and
$c(n_2,m_2)$ can be rewritten
\begin{eqnarray*}
h_G(homogeneous) &=&-\frac J{4m}(m_1+m_2+m_{12})+\frac{\gamma ^{*}}{4m}%
N(N-1) \\
&=&-\frac J{4m}(m_1+m_2+m_{12})+\frac{\gamma ^{*}}{4m}\left[
n_1(n_1-1)+n_2(n_2-1)+2n_1n_2\right] \\
&>&-\frac J{4m}(m_1+m_2)+\frac{\gamma ^{*}}{4m}\left[
n_1(n_1-1)+n_2(n_2-1)\right] \\
&=&h_G(diverse)
\end{eqnarray*}
which is in turn equivalent to the normalized value of the outlink
density or inter-community links density,
\begin{eqnarray*}
\frac{Jm_{12}}{2n_1n_2} &=&\gamma ^{*}, \\
&&
\end{eqnarray*}
since $m=m_1+m_2+m_{12}$ and $N=n_1+n_2.$

The proofs of LDP over sparse Erdos Renyi random graphs $ER(N,p)$
and scale-free Chung-Lu graphs $CL_\infty (N,\varpi )$ with expected
degrees sequence $\varpi $ will be given for arbitrary $f$ in
(\ref{h}) since they
do not depend on the form of the function $f.$ To begin the proof for $%
ER(N,p),$ based on the Azuma-Hoefding inequalities, we define, as before, $%
H(G)\equiv \min h_G=h_G(\vec{s}_0)$ and label the optimum
$\vec{s}_0.$ We need a technical estimate required in the
application of theorem 1.

\textbf{Lemma 3:} There is a constsnt $c>0$ independent of $N$ such
that for any two graphs $G(N,m^{\prime }),$ $G^{\prime
}(N,m^{^{\prime \prime }}=m^{\prime }+1)$ in $random$ ensemble
$\Omega (N)$ which differs in exactly one edge,
\[
|H(G)-H(G^{\prime })|<c
\]

Proof: For a fixed state $\vec{s}=(s_1,...,s_N),$ since the second
term $f$ cancel,
\begin{eqnarray*}
\frac{4m^{\prime }}J|h_G(\vec{s})-h_{G^{\prime }}(\vec{s})|
&=&\left| \sum_{i\neq j}\left[ A_{ij}-A_{ij}^{^{\prime }}\right]
\delta
(s_i,s_j)\right| \\
&=&\left| \left[ A_{i_0j_0}-A_{i_0j_0}^{^{\prime }}\right] \delta
(s_{i_0},s_{j_0})\right| \\
&\leq &2,
\end{eqnarray*}
implying that $|h_G(\vec{s})-h_{G^{\prime }}(\vec{s})|\leq \frac
J{4m}.$ A standard argument based on assumimg wlog $H(G)\leq
H(G^{\prime }),$ that is, $h_G(\vec{s}_0)=\min h_G=H(G),$
$h_{G^{\prime }}(\vec{s}_0^{^{\prime
}})=\min h_{G^{\prime }}=H(G^{\prime }),$ and $h_{G^{\prime }}(\vec{s}%
_0)\geq \min h_{G^{\prime }}\equiv H(G^{\prime }),$ so $h_G(\vec{s}%
_0)=H(G)\leq H(G^{\prime })\leq h_{G^{\prime }}(\vec{s}_0),$ implies
the result,
\begin{eqnarray*}
|H(G)-H(G^{\prime })| &<&|h_G(\vec{s}_0)-h_{G^{\prime }}(\vec{s}_0)| \\
&\leq &\frac J2.
\end{eqnarray*}

The rest of the proof is exactly the same as in previous sections -
using the Divide and Conquer technique based on the non-sparse
property which is valid for $ER[p,N]$ and $CL_\infty [N,B,\beta ]$
and theorem 1 together with this lemma, we prove:

\textbf{Theorem 10} \emph{In the non-sparse $ER[N,p]$ and }$CL_\infty
[N,B,\beta ]$\emph{\ ensembles $\Omega $, the optimal Q-Potts functional $%
H(G)$ satisfies the concentration inequalities respectively: there exists }$%
\mu >0$ independent of $N$ such that for real $t>0,$

\[
P(|H(G)-<H>_\Omega |>t)\leq 2\exp {(-{\mu ^2N^2}-{t^2})}+2\exp {(-{%
t^2N^2(p-\mu )^2})}
\]

\[
P(|H(G)-<H>_\Omega |>t)\leq 2\exp {(-{\frac{\mu ^2N^{2\beta }}{2\lambda ^2}}-%
{\frac{t^2}{2\sigma ^2}})}+{2\exp }\left( -\frac{t^2(B/2-\mu )^2N^{2\beta }}%
2\right)
\]
with $\lambda =1/2$ and $\sigma =J/2.$

Using Lemma 3 and Theorem 3, we prove:

\textbf{Theorem 11} \emph{In $ER[Np,K]$ and }$CL[N,\varpi
,K]$\emph{\ ensembles $\Omega $ with uniform upper bound $K$ on the
degree which is independent of }$N$\emph{, the optimal Q-Potts
functional $H(G)$ satisfies the concentration inequality: for any
real }$t>0,$
\[
\mathscr{P}(|H(G)-<H>_\Omega|>t)\leq 2\exp {(-{\frac{Nt^2}{8K^2}}).}
\]

\section{Conclusion}
In this paper, we derive LDP type inequalities for the optimal values of community detecting functions on random networks to show its concentration and stability. There is no concentration for the most general case. We prove the concentration of the general community detecting function on K-bound sparse network and an even sharper concentration when the network ensemble is generated by a given network and a small perturbation. Then we examine several specific cases.
The three specific community detecting functions we considered are: modularity, graph bipartioning and q-potts community structure. The specific network types we considered are ER and CL networks, and each of them with sparse(K-bound) and non-sparse cases. We prove concentration in these cases, that means in these cases the community detecting functions are stable especially when the system size is large enough.

%\begin{acknowledgements}
This work was supported in part by the Army Research Laboratory
under Cooperative Agreement Number W911NF-09-2-0053, and the Army
Research Office Grant No. W911NF-09-1-0254. The views and conclusions
contained in this document are those of the authors and should not
be interpreted as representing the official policies, either
expressed or implied, of the Army Research Laboratory or the U.S.
Government.
%\end{acknowledgements}

\end{document}